\documentclass[10pt]{amsart}
\usepackage[all]{xy}
\usepackage{setspace}
\newlength{\baseunit}               
\newcommand{\bpf}{\noindent {\em Proof.  }}
\newcommand{\epf}{\qed \vspace{+10pt}}
\newtheorem{tm}{Theorem}
\newtheorem{pr}[tm]{Proposition}
\newtheorem{lm}[tm]{Lemma}
\newtheorem{co}[tm]{Corollary}

\newcommand{\proj}{\mathbb P}

\newcommand{\oh}{{\mathcal{O}}}

\newcommand{\Spec}{\operatorname{Spec}}

\begin{document}
\pagestyle{headings}
\title{Relative Canonical sheaves of a Family of Curves}
\author{Jongmin Lee}

\begin{abstract} In this paper we show that the relative canonical
sheaf of a relatively minimal fibration of curves over a curve is
semi-ample ; in fact, its m-tensored product is base point free for any m
$\ge 2$.  We use Koszul cohomology with it to prove that the
relative canonical ring of the fibration is generated in degree up
to five.
\end{abstract}

\maketitle

\section{Introduction}

Let $C$ be a nonhyperelliptic projective nonsingular curve of
genus $g$ over $\mathbb{C}$. Then the theorem of Noether says that
$$H^0(C,\omega_C) \otimes H^0(C,\omega_C) \to H^0(C,\omega_C^{\otimes 2})$$
is surjective. Hence, the canonical ring $R_C=\bigoplus_{m \ge 0}
H^0(C,\omega_C^{\otimes m})$ of $C$ is generated by elements of
degree $\le 1$. If $C$ is hyperelliptic, direct computations show
that $R_C$ is generated by elements of degree $\le 2$. Recently,
K. Konno(\cite{konno}) proved the following :

Let $F$ be a fiber in a relatively minimal fibration of curves of
genus $g \ge 2$ over $\mathbb{C}$. Then the canonical ring
$R(F,K_F)$ is generated in degrees 1,2 and 3 except when $F$ is a
multiple fiber which contains a $(-1)$-elliptic cycle $E$ such
that $E \subset Bs|K_F|$. In the exceptional case, it needs one
more generator in degree 4.

In particular, when $f \colon X \rightarrow S$ is a relatively
minimal nonsingular surface over a nonsingular curve $S$ over
$\mathbb{C}$, the relative canonical ring $\mathcal{R}(f) =
\bigoplus_{m \ge 0} f_*(\omega_{X/S}^{\otimes m})$ is generated in
degree $\le 4$. In this paper we study the relative canonical ring
$\mathcal{R}(f)$ when $X,S$ and $f$ are not necessarily over
$\mathbb{C}$. We allow $S$ to be the spectrum of any discrete
valuation ring. The main ingredient of the study consists of
showing the global generatedness of $\omega_{X/S}^{\otimes m}$'s
for $m \ge 2$, which can be shown easily over $\mathbb{C}$ using
the proposition 2.3 of F. Catanese and M. Franciosi in \cite{cf} :

Let $C$ be a curve lying on a smooth algebraic surface $X$ over an
algebraically closed field and let $D$ be a divisor on $C$. A
point $x \in C$ is not a base point for $|D|$ if for every
subcurve $Y$ of $C$ we have $D \cdot Y \ge 2 p_a(Y)$.

We also use Koszul cohomology(\cite{konno} or \cite{green}) to
conclude that $\mathcal{R}(f)$ is generated in degree $\le 5$.

\section{Statement of the Problem}

\subsection{Settings}

Let $X \to S$ be a flat proper morphism satisfying the following.

\begin{enumerate}
\item $S$ is the spectrum of a discrete valuation ring $R$ with
uniformizer $t$ and residue field $k=R/(t)$.
\item The generic fiber of $X \to S$ is a geometrically
irreducible curve of genus $g \geq 2$. \item $X$ is a two
dimensional regular scheme.
\item $X$ is relatively minimal over $S$, i.e.,there are no
$(-1)$-curves contained in the special fiber of $X \to S$.
\end{enumerate}

\noindent To understand the last condition we need to define what
$(-1)$-curves are. Let $C \subset X$ be an integral curve in the
special fiber with $n=[\Gamma(C,\oh_C):k]$. If $C$ is a rational
curve over $\Gamma(C,\oh_C)$ with $C^2=-dn$, we call $C$ a
$(-d)$-curve on $X$. Here the intersection number is over $k$.

\noindent We will use the notation $X_0 = X \times_S $Spec$(k)$
for the special fibre.

\subsection{Object}
We are going to prove that $\omega_{X/S}^{\otimes m}$ is globally
generated for $m \ge 2$.

\subsection{Observation}
Let $R^{sh}$ be the strict Henselization of $R$,
$S^{sh}=\Spec(R^{sh})$ and $X^{sh}=X \times_S S^{sh}$. Then we
have $H^0(X^{sh},\omega_{X^{sh}/S^{sh}}^{\otimes n}) \cong
H^0(X,\omega_{X/S}^{\otimes n})\otimes_R R^{sh}$ and $X^{sh}$ is
regular as well. Therefore we can assume that $k$ is separably
closed. One good thing with this assumption is that any
$(-2)$-curve over a separably closed field $k$ is isomorphic to
the projective line over $k$(it can be shown by a little
adaptation of Lemma 3.4 in \cite{lich} together with the existence
of a $k$-rational point on the curve). Also note that any finite
extension of a separably closed field is separably closed.
Therefore any $(-2)$-curve $C$ on $X$ can be assumed to be
isomorphic to the projective line over $\Gamma(C,\oh_C)$.

\section{Lemmas}

\begin{lm} \label{conic}
Let $C$ be an integral projective curve over a field $k$. If
$deg(\omega_C)<0$, then $C\subset\proj_{k'}^2$ is an irreducible
conic where $k'=\Gamma(C,\oh_C)$.
\end{lm}

\bpf We may assume that $k'=k$ and that $\deg(\omega_C)=-2$. It is
because $\deg(\omega_C)$ is $\geq -2$ and even. Then
$\chi(\omega_C^{-1})=\deg(\omega_C^{-1})+\chi(\oh_C)=2+1=3$ and
hence $h^0(C,\omega_C^{-1})\ge3$. Any nonzero section $s$ of
$\omega_C^{-1}$ has a zero divisor which is an effective Cartier
divisor $D$ of degree 2. Thus for reason of dimension in the
sequence
$$0 \longrightarrow \oh_C \longrightarrow \omega_C^{-1} \longrightarrow \oh_D \longrightarrow 0$$
the map $H^0(C,\omega_C^{-1}) \to H^0(\oh_D)$ is surjective. Hence
$\omega_C^{-1}$ is generated by global sections. Then
$H^0(C,\omega_C^{-1})$ defines an isomorphism between $C$ and a
conic in $\proj^2$. \epf

\begin{lm} \label{van0}
Let $C$ be an integral projective Gorenstein curve over $k$. Let
$\mathcal{L}$ be an invertible $\oh_C$-module on $C$ with
$\deg(\mathcal{L}) > \deg(\omega_C)$. Then $H^1(C,\mathcal{L})=0$
and one of the following two possibilities occurs:
$H^0(C,\mathcal{L}) \neq 0$ or $C \cong \proj_{k'}^1$ for some
finite extension $k'$ and $\mathcal{L}=\oh_{\proj^1}(-1)$.
\end{lm}

\bpf The vanishing of $H^1$ follows from Serre duality. If
$h^0(C,\mathcal{L})=0$ then $\chi(\mathcal{L})=0$. Hence $
\deg(\mathcal{L})=\frac12 \deg(\omega_C)$. Then $
\deg(\mathcal{L})
> \deg(\omega_C)$ implies that $ \deg(\omega_C) < 0$. By the previous
lemma $C \subset \proj_{k'}^2$ is a conic. The existence of the
degree $-1$ invertible sheaf $\mathcal{L}$ implies that $C$ is a
smooth conic with a rational point. This follows by considering
the linear series $|\mathcal{L}^{-1}|$. \epf

\begin{lm} \label{gg}
Let $C$ be an integral Gorenstein curve over $k$ and let
$k'=\Gamma(C,\oh_C)$. Let $\mathcal{L}$ be an invertible
$\oh_C$-module with $$\deg(\mathcal{L}) \geq \deg(\omega_C) +
2[k':k].$$ Then $\mathcal{L}$ is globally generated.
\end{lm}

\bpf Let $B \subset C$ be the base locus of $|\mathcal{L}|$. Since
$H^0(C,\mathcal{L}) \neq 0$ by Riemann-Roch we see that $B \neq
C$. The exact sequence $$0 \longrightarrow I_B \otimes \mathcal{L}
\longrightarrow \mathcal{L} \longrightarrow \oh_B \longrightarrow
0$$ shows that $H^1(C,I_B \otimes \mathcal{L}) \neq 0$ if $B \neq
\emptyset$. By duality, this implies that there is a nonzero map
$I_B \otimes \mathcal{L} \to \omega_C$. In particular
$$h^0(\mathcal{L})=h^0(I_B \otimes \mathcal{L}) \leq
h^0(\omega_C).$$ This is impossible by our assumption. \epf

\noindent The following is a slight generalizion of Lemma 2.1 in
\cite{cf}.

\begin{lm} \label{van}
Let $C$ be a curve lying on a smooth surface $X$(i.e. $C$ is an
effective divisor on $X$) and let $D$ be a divisor on $C$. Then
$H^1(C,\oh_C(D))=0$ if for all subcurves $B < C$, $D \cdot B >
deg(\omega_B)$.
\end{lm}

\bpf If $H^1(C,\oh_C(D)) \cong H^0(\oh_C(K_C-D)) \neq 0$ there is
a nonzero section $\sigma \in H^0(\oh_C(K_C-D))$. Let $Z$ be the
maximal curve in $C$ on which $\sigma$ vanishes identically and
let $Y=C-Z$. Then by the exact sequence $$0 \longrightarrow
\oh_Y(K_C-D-Z) \longrightarrow \oh_C(K_C-D) \longrightarrow
\oh_Z(K_C-D) \longrightarrow 0,$$ upon dividing $\sigma$ by a
section $\zeta$ with $ $div$(\zeta)=Z$, we obtain
$\frac{\sigma}{\zeta}=\sigma'$ a section of
$H^0(Y,\oh_Y(K_C-D-Z))$ vanishing on a finite set. Hence, we have
$(K_C-D) \cdot Y \geq Z \cdot Y$. This inequality, since by
adjunction $\oh_C(K_C)=\oh_C(K_X+C)$,is equivalent to $D \cdot Y
\leq \deg(\omega_Y)$, a contradiction. \epf

\section{Configurations of (-2)-chains and Fundamental Cycles}
Let $\sum_{i=1}^s r_i C_i$ be a connected chain of $(-2)$-curves
in $X_0$. Here $r_i$ is the multiplicity of $C_i$ in $X_0$. Let us
use the notation $k_i=\Gamma(C_i,\oh_{C_i})$.

\begin{lm} \label{strict}
Let $n_i=[k_i:k]$. If $C_i \cdot C_j > 0$, then we have $C_i \cdot
C_j = \max(n_i,n_j)$.
\end{lm}
\bpf Let us suppose that $n_i \geq n_j$. If the claim is false,
then we have $C_i \cdot C_j \geq 2n_i$. It implies that
$(C_i+C_j)^2=-2n_i+2C_i \cdot C_j-2n_j \geq 2n_i-2n_j$. We also
have $(C_i+C_j)^2 < 0$ since $X_0$ has genus $\geq 2$ by
assumption so there is a curve in $X_0$ which is not a
$(-2)$-curve. Hence we have a contradiction $n_i < n_j$. \epf
\begin{lm}
If $C_i \cdot C_j > 0$ with $n_i \geq n_j$, then $n_i=n_j,2n_j$ or
$3n_j$.
\end{lm}
\bpf It is clear from $0 > (C_i+2C_j)^2=-2n_i+4C_i \cdot C_j-8n_j
= 2n_i-8n_j$. The first inequality is strict by the same reason in
Lemma \ref{strict} \epf

Let's denote $<C_i,C_j>=-2\frac{C_i \cdot C_j}{C_i^2}$. Then
$(<C_i,C_j>)_{i,j=1,...,s}$ has the same properties as Cartan
matrices associated to Lie algebras(refer chapter 4 in \cite{ja}).
Therefore we have the following list of possible dual graphs for
connected $(-2)$-chains. The number attached below each $C_i$ is
$[k_i:k]$.


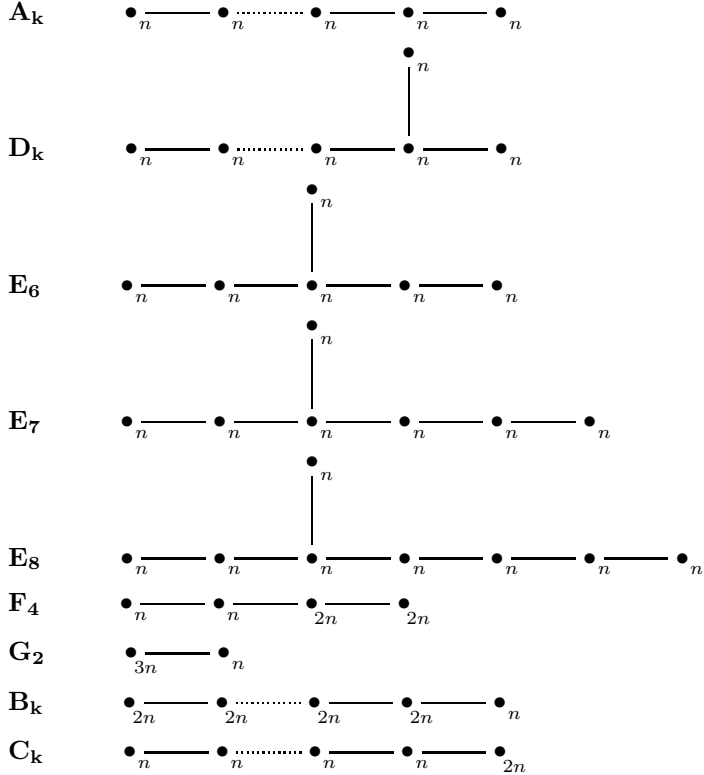
\xymatrix{ {\bf A_k} \ar@{}[r] & {\bullet} \ar@{-}[r]_<{n} &
{\bullet} \ar@{.}[r]_<{n} & {\bullet} \ar@{-}[r]_<{n} & {\bullet}
\ar@{-}[r]_<{n} & {\bullet} \ar@{}[r]_<{n} & }

\xymatrix{ & & & & {\bullet} &
\\
{\bf D_k} \ar@{}[r] & {\bullet} \ar@{-}[r]_<{n} & {\bullet}
\ar@{.}[r]_<{n} & {\bullet} \ar@{-}[r]_<{n} & {\bullet}
\ar@{-}[r]_<{n} \ar@{-}[u]_>{n} & {\bullet} \ar@{}[r]_<{n} & }

\xymatrix{ & & & {\bullet} &
\\
{\bf E_6} \ar@{}[r] & {\bullet} \ar@{-}[r]_<{n} & {\bullet}
\ar@{-}[r]_<{n} & {\bullet} \ar@{-}[r]_<{n} \ar@{-}[u]_>{n} &
{\bullet} \ar@{-}[r]_<{n} & {\bullet} \ar@{}[r]_<{n} & }
\xymatrix{ & & & {\bullet} & \\
{\bf E_7}  \ar@{}[r] & {\bullet} \ar@{-}[r]_<{n} & {\bullet}
\ar@{-}[r]_<{n} & {\bullet} \ar@{-}[r]_<{n} \ar@{-}[u]_>{n} &
{\bullet} \ar@{-}[r]_<{n} & {\bullet} \ar@{-}[r]_<{n} & {\bullet}
\ar@{}[r]_<{n} & }
\xymatrix{ & & & {\bullet} \\
{\bf E_8} \ar@{}[r] & {\bullet} \ar@{-}[r]_<{n} & {\bullet}
\ar@{-}[r]_<{n} & {\bullet} \ar@{-}[r]_<{n} \ar@{-}[u]_>{n} &
{\bullet} \ar@{-}[r]_<{n} & {\bullet} \ar@{-}[r]_<{n} & {\bullet}
\ar@{-}[r]_<{n} & {\bullet} \ar@{}[r]_<{n} & }

\xymatrix{ {\bf F_4} \ar@{}[r] & {\bullet} \ar@{-}[r]_<{n} &
{\bullet} \ar@{-}[r]_<{n} & {\bullet} \ar@{-}[r]_<{2n} & {\bullet}
\ar@{}[r]_<{2n} & }

\xymatrix{ {\bf G_2} \ar@{}[r] & {\bullet} \ar@{-}[r]_<{3n} &
{\bullet} \ar@{}[r]_<{n} & }

\xymatrix{ {\bf B_k} \ar@{}[r] & {\bullet} \ar@{-}[r]_<{2n} &
{\bullet} \ar@{.}[r]_<{2n} & {\bullet} \ar@{-}[r]_<{2n} &
{\bullet} \ar@{-}[r]_<{2n} & {\bullet} \ar@{}[r]_<{n} & }

\xymatrix{ {\bf C_k} \ar@{}[r] & {\bullet} \ar@{-}[r]_<{n} &
{\bullet} \ar@{.}[r]_<{n} & {\bullet} \ar@{-}[r]_<{n} & {\bullet}
\ar@{-}[r]_<{n} & {\bullet} \ar@{}[r]_<{2n} &
}


\vskip 5mm

Notice here that all types can be appeared in any characteristic
except that $F_4,B_k$ and $C_k$ are only in char 2 and $G_2$ is
only in char 3.

Consider all nonzero effective divisors $Z$ of the form
$\sum_{i=1}^s a_i C_i$ which are less than or equal to
$\sum_{i=1}^s r_i C_i$ and satisfies that $Z \cdot C_i \leq 0$ for
all $i=1,...,s$. We will call the minimal elements among them the
\emph{fundamental cycles}(refer \cite{bpv}). The following dual
graphs are the fundamental cycles associated to the connected
$(-2)$-chains we listed above. The number above each $C_i$ is the
multiplicity of $C_i$ in the fundamental cycle.

\vskip 1cm

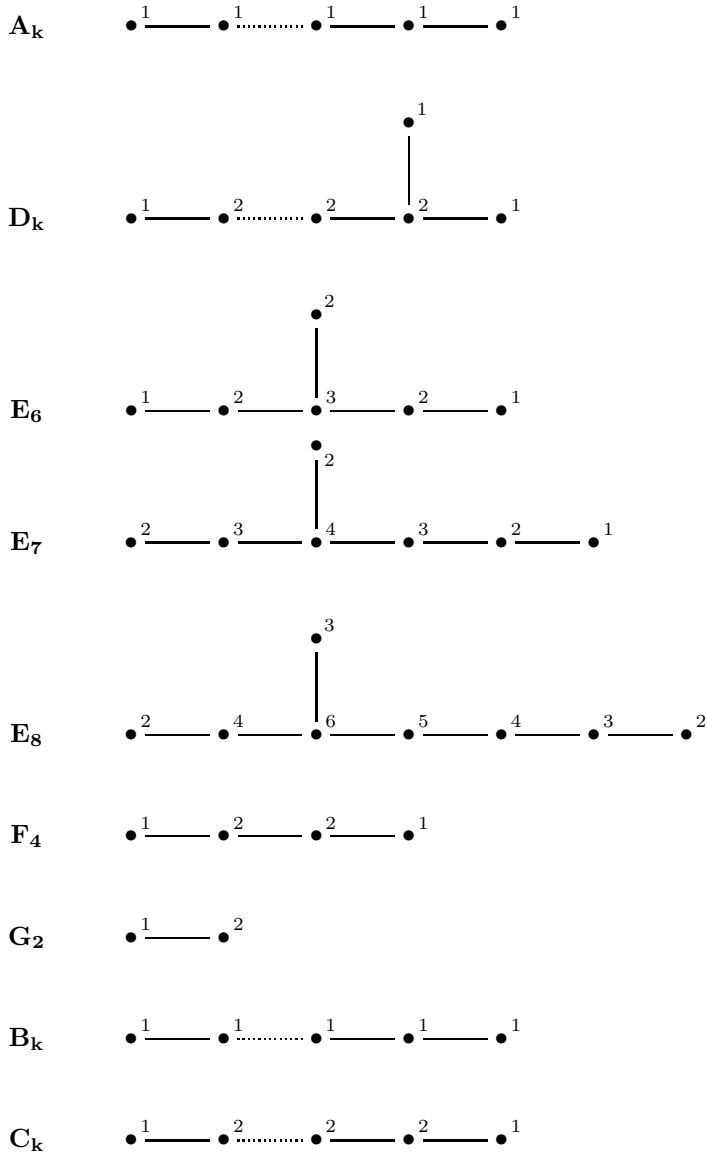
\xymatrix{ {\bf A_k} & {\bullet} \ar@{-}[r]^<{1} & {\bullet}
\ar@{.}[r]^<{1} & {\bullet} \ar@{-}[r]^<{1} & {\bullet}
\ar@{-}[r]^<{1} & {\bullet} \ar@{}[r]^<{1} &
\\
& & & & {\bullet} \ar@{}[u]_<{1}
\\
{\bf D_k} & {\bullet} \ar@{-}[r]^<{1} & {\bullet} \ar@{.}[r]^<{2}
& {\bullet} \ar@{-}[r]^<{2} & {\bullet} \ar@{-}[r]^<{2} \ar@{-}[u]
& {\bullet} \ar@{}[r]^<{1} &
\\
& & & {\bullet} \ar@{}[u]_<{2}
\\
{\bf E_6} & {\bullet} \ar@{-}[r]^<{1} & {\bullet} \ar@{-}[r]^<{2}
& {\bullet} \ar@{-}[r]^<{3} \ar@{-}[u] & {\bullet} \ar@{-}[r]^<{2}
& {\bullet} \ar@{}[r]^<{1} & }
\xymatrix{ & & & {\bullet} \ar@{}[d]^<{2} \\
{\bf E_7} & {\bullet} \ar@{-}[r]^<{2} & {\bullet} \ar@{-}[r]^<{3}
& {\bullet} \ar@{-}[r]^<{4} \ar@{-}[u] & {\bullet} \ar@{-}[r]^<{3}
& {\bullet} \ar@{-}[r]^<{2} & {\bullet} \ar@{}[r]^<{1} & \\ & & &
{\bullet} \ar@{}[u]_<{3} \\ {\bf E_8} & {\bullet} \ar@{-}[r]^<{2}
& {\bullet} \ar@{-}[r]^<{4} & {\bullet} \ar@{-}[r]^<{6} \ar@{-}[u]
& {\bullet} \ar@{-}[r]^<{5} & {\bullet} \ar@{-}[r]^<{4} &
{\bullet}
\ar@{-}[r]^<{3} & {\bullet} \ar@{}[r]^<{2} & \\
{\bf F_4} & {\bullet} \ar@{-}[r]^<{1} & {\bullet} \ar@{-}[r]^<{2}
& {\bullet} \ar@{-}[r]^<{2} & {\bullet} \ar@{}[r]^<{1} &
\\
{\bf G_2} & {\bullet} \ar@{-}[r]^<{1} & {\bullet} \ar@{}[r]^<{2} &
\\
{\bf B_k} & {\bullet} \ar@{-}[r]^<{1} & {\bullet} \ar@{.}[r]^<{1}
& {\bullet} \ar@{-}[r]^<{1} & {\bullet} \ar@{-}[r]^<{1} &
{\bullet} \ar@{}[r]^<{1} &
\\
{\bf C_k} & {\bullet} \ar@{-}[r]^<{1} & {\bullet} \ar@{.}[r]^<{2}
& {\bullet} \ar@{-}[r]^<{2} & {\bullet} \ar@{-}[r]^<{2} &
{\bullet} \ar@{}[r]^<{1} & }

\vskip 5mm

\noindent Notice here that each fundamental cycle $Z$ satisfies
$Z^2=-2n$, where $n$ is the minimum of $n_i$'s.

\section{Base Point Freeness}

\begin{tm} \label{main}
$\omega_{X/S}^{\otimes m}$ is globally generated for $m \geq 2$.
\end{tm} \bpf

We only need to show that $\omega_{X/S}^{\otimes m}$ is globally
generated on $X_0$ since any closed point is in $X_0$ by our
assupmtion. Let us consider the following exact sequence
$$0 \longrightarrow \omega_{X/S}^{\otimes m}(-X_0) \longrightarrow
\omega_{X/S}^{\otimes m} \longrightarrow \omega_{X/S}^{\otimes
m}|_{X_0} \longrightarrow 0.$$ Using Lemma \ref{van}, we have
$H^1(X_0,\omega_{X/S}^{\otimes m}(-X_0)|_{X_0})=0$ for any $m \geq
2$ . Then Nakayama lemma implies that $H^1(X,\omega_{X/S}^{\otimes
m}(-X_0))=0$. Hence, we get a surjective restriction map
$H^0(X,\omega_{X/S}^{\otimes m}) \to H^0(X_0,\omega_{X/S}^{\otimes
m}|_{X_0})$. Therefore, it is sufficient to show that
$\omega_{X/S}^{\otimes m}|_{X_0}$ is globally generated. We are
going to prove it by induction. Let $Y \subset X_0$ be a divisor.
Suppose that for all proper subdivisors of $Y' \subset Y$ the
restriction $\omega_{X/S}^{\otimes m}|_{Y'}$ is globally
generated. Let $C \subset Y$ be an irreducible component of $Y$.
Set $Y'=Y-C$ and consider the exact sequence $$0 \longrightarrow
\omega_{X/S}^{\otimes m}(-Y')|_C \longrightarrow
\omega_{X/S}^{\otimes m}|_{Y} \longrightarrow
\omega_{X/S}^{\otimes m}|_{Y'} \longrightarrow 0$$ of coherent
sheaves. We compute
\begin{eqnarray}
\nonumber \deg_C(\omega_{X/S}^{\otimes m}(-Y')|_C)=mK_{X/S} \cdot
C - Y \cdot C + C^2 \\
\nonumber = \deg_C(\omega_C) + (m-1)K_{X/S} \cdot C - Y \cdot C.
\end{eqnarray}
With this in mind we choose $C$ as follows: if $Y^2<0$ then we
choose $C$ such that $C \cdot Y < 0$(such $C$ exists since
otherwise $Y^2 \geq 0$), if $Y^2=0$ then we choose $C$ such that
$K_{X/S} \cdot C>0$(if there is no such $C$, we have $K_{X/S}
\cdot Y=0$;hence $K_{X/S} \cdot X_0=0$). In either case, by Lemma
\ref{van0}, we'll have $H^1(C,\omega_{X/S}^{\otimes
m}(-Y')|_C)=0$. Therefore the restriction map
$H^0(Y,\omega_{X/S}^{\otimes m}|_Y) \to
H^0(Y',\omega_{X/S}^{\otimes m}|_{Y'})$ is surjective. And we
conclude that $\omega_{X/S}^{\otimes m}|_Y$ is globally generated
at all points of $Y'$. If $C$ occurs in $Y$ with multiplicity
$\geq 2$ then the fact that $\omega_{X/S}^{\otimes m}|_{Y'}$ is
globally generated implies the result for $\omega_{X/S}^{\otimes
m}|_Y$. Hence we can assume that $C$ occurs with multiplicity 1 in
$Y$. There are several cases to consider. From now on we will use
the notation $k'=\Gamma(C,\oh_C)$.

\noindent{\bf Case 1.} $K_{X/S} \cdot C=0$ and $C \cap Y' \neq
\emptyset$. \vskip 0.1cm \noindent Let $s'$ be a section of
$\omega_{X/S}^{\otimes m}|_{Y'}$ which does not vanish identically
on $C \cap Y'$. Let $s$ be a section of $\omega_{X/S}^{\otimes
m}|_Y$ which lifts $s'$. Since $\omega_{X/S}^{\otimes m}|_C$ has
degree $0$ we conclude that the restriction $s|_C$ of $s$ is a
nonvanishing section on $C$. This proves that
$\omega_{X/S}^{\otimes m}|_Y$ is globally generated in this case.

\noindent{\bf Case 2.} $K_{X/S} \cdot C=0$ and $C \cap Y' =
\emptyset$. \vskip 0.1cm \noindent We have $Y=C \amalg Y'$ and
$\omega_{X/S}^{\otimes m}(-Y')|_C=\omega_{X/S}^{\otimes m}|_C$.
This is a sheaf of degree $ $deg$_C(\omega_{X/S}^{\otimes
m}|_C)=mK_{X/S} \cdot C = 0$. On the other hand we also have $
$deg$(\omega_C)=K_{X/S} \cdot C + C^2 = C^2 < 0$(because if
$C^2=0$ then $X_0=rC$ and $K_{X/S} \cdot X_0 = rK_{X/S} \cdot C =
0$, a contradiction to the assumption $g \geq 2$). Thus we
conclude by Lemma \ref{conic} that $C \subset \proj_{k'}^2$ is a
conic. Since the Picard group of an irreducible conic is
$\mathbb{Z}$, $ \omega_{X/S}|_C \cong \oh_C$.

\noindent{\bf Case 3.} $K_{X/S} \cdot C > 0$ and $m \ge 3$. \vskip
0.1cm \noindent Due to the way we chose $C$, we have $Y \cdot C
\leq 0$. Then $\omega_{X/S}^{\otimes m}(-Y')|_C$ satisfies the
degree condition in Lemma \ref{gg} hence it is globally generated.
It implies that $\omega_{X/S}^{\otimes m}|_Y$ is globally
generated on $C$.

\noindent{\bf Case 4.} $K_{X/S} \cdot C > 0$, $Y \cdot C < 0$ and
$m=2$. \vskip 0.1cm \noindent Lemma \ref{gg} implies again that
$\omega_{X/S}^{\otimes 2}(-Y')|_C$ is globally generated.

\noindent{\bf Case 5.} $K_{X/S} \cdot C
>0$, $Y \cdot C = 0$ and $m=2$; hence we have $Y^2=0$. \vskip 0.1cm
\noindent If there are more than one such curve $C$, then
$\omega_{X/S}^{\otimes 2}|_Y$ is globally generated by choosing
$C$ alternately. Hence we can assume that we have only one $C$
with $K_{X/S} \cdot C >0$. If $K_{X/S} \cdot C \geq 2[k':k]$, we
are done by Lemma \ref{gg}. If not then $K_{X/S} \cdot C = [k':k]$
and we see that $ $deg$(\omega_C)=K_{X/S} \cdot C + C^2$ is either
(a) $0$ or (b) $-2[k':k]$.

\noindent{\bf Case 5(a).} We see that $0=Y \cdot C=C^2+C \cdot Y'$
implies that the scheme theoretic intersection $C \cap Y'$ is a
single point $p$ with residue field $k'$. By Riemann-Roch we have
$ $dim$_{k'}H^0(C,\omega_{X/S}^{\otimes 2}|_C)=2$. On the other
hand, the restriction map $\omega_{X/S}^{\otimes 2}|_Y \to
\omega_{X/S}^{\otimes 2}|_C$ gives us a section $t$ of
$\omega_{X/S}^{\otimes 2}|_C$ which does not vanish at $p$
(namely, lift a section $s'$ of $\omega_{X/S}^{\otimes 2}|_{Y'}$
which does not vanish at $p$ to a section $s$ of
$\omega_{X/S}^{\otimes 2}|_Y$ and then restrict this to get $t$ on
$C$). In addition $\omega_{X/S}^{\otimes 2}(-Y')|_C$ has degree
$[k':k]$ and thus by Riemann-Roch has a nonzero section $t'$. We
can think of $t'$ as a section of $\omega_{X/S}^{\otimes 2}|_C$
because of the short exact sequence
$$0 \longrightarrow \omega_{X/S}^{\otimes 2}(-Y')|_C \longrightarrow
\omega_{X/S}^{\otimes 2}|_C \longrightarrow \oh_p \longrightarrow
0.$$ Clearly $t'$ vanishes at $p$ and hence ${t,t'}$ is a
$k'$-basis of $H^0(C,\omega_{X/S}^{\otimes 2}|_C)$. By Lemma
\ref{gg} applied to $\omega_{X/S}^{\otimes 2}|_C$ on $C$ we see
that $t$ and $t'$ do not have a common zero on $C$. We conclude
that the image of $H^0(Y,\omega_{X/S}^{\otimes 2}|_Y) \to
H^0(C,\omega_{X/S}^{\otimes 2}|_C)$ generates
$\omega_{X/S}^{\otimes 2}|_C$.

\noindent{\bf Case 5(b).} In this case we have $C^2=-3[k':k]$ and
hence the conic $C \subset \proj^2_{k'}$ has an invertible sheaf
of degree 1. We conclude that $C \cong \proj^1_{k'}$. Therefore we
have $Y=C+\sum_{i \in I} a_i C_i$ where $C^2=-3[k':k], K_{X/S}
\cdot C = [k':k]$ and $C_i$'s are $(-2)$-curves(Notice that $\deg
\omega_Y = Y^2 + K \cdot Y = [k':k]$. Hence $[k':k]$ is even so
this case should have char($k$)=$2$. If you don't have any
interest in positive characteristics you can skip this case). By
induction hypothesis we know that $\omega_{X/S}^{\otimes 2}|_{Y'}$
is globally generated and it implies that $\omega_{X/S}^{\otimes
2}|_Y$ is globally generated on $Y'$. Therefore we only need to
show that $\omega_{X/S}^{\otimes 2}|_Y$ is globally generated on
$C-\cup_{i \in I} C_i$. We are going to find appropriate sections
by gluing.

Since we have $C \cdot \sum_{i \in I} a_i C_i = 3[k':k]$, we have
at most three intersection points between $C$ and $\cup_{i \in I}
C_i$ and also at most three connected components of $\cup_{i \in
I} C_i$. Let us write $D=\sum_{i \in I} a_i C_i$. Let $\alpha$ be
the number of connected components of $D$ and $\beta$ that of the
set theoretic intersection points $C \cap D$.

(I) Case $\alpha=3$.

We have $\beta=3$. Let $x_1$,$x_2$ and $x_3$ be the three
intersection points. Let $D_1$,$D_2$ and $D_3$ be the connected
components where $D_i$ meets $C$ at $x_i$. Take the zero sections
for $D_1$ and $D_2$ and take a section $\sigma$ for $D_3$ with
$\sigma(x_3) \not= 0$(i.e. take the restriction of a section of
$\omega_{X/S}^{\otimes 2}|_{Y'}$ not vanishing at $x$). We can
find a section $\tau$ of $\omega_{X/S}^{\otimes 2}|_C$ with
$\tau(x_1)=\tau(x_2)=0$ and $\tau(x_3)=\sigma(x_3)$. Then these
sections give a section of $\omega_{X/S}^{\otimes 2}|_Y$ and it
does not vanish on C except at $x_1$ and $x_2$.

(II) Case $\alpha=2$.

We have $\beta=2$ or $3$. Let $D_1$ and $D_2$ be the connected
components where $C \cap D_2$ is a rational point $x$ on $C$
scheme theoretically. Take the zero section for $D_1$ and take a
section $\sigma$ for $D_3$ with $\sigma(x) \not= 0$. Then we can
find a section $\tau$ of $\omega_{X/S}^{\otimes 2}|_C$ such that
$\tau|_{C \cap D_1 }\equiv0$ and $\tau(x)=\sigma(x)$. Then these
sections give a section of $\omega_{X/S}^{\otimes 2}|_Y$ and it
vanishes only on $C \cap D_1$.

(III) Case $\alpha=1$.

Let $m=[k':k]$ and $Z$ the fundamental cycle of $D$. Then we have
\begin{equation*} \label{eq*}
\tag{$*$} 0 \geq (C+Z)^2 = C^2+Z^2+2C \cdot Z = -3m-2n+2C \cdot Z.
\end{equation*}
Here $n$ is the smallest among $[k(C_i):k]$'s. If $\beta=3$, we
have $m \ge n$ and $C \cdot Z=3m$ which contradict (\ref{eq*}).
Therefore $\beta=1$ or $2$.

(i) Case $\beta=2$.

We have $m \ge n$ and $C \cdot Z=3m$ or $2m$. We can assume $C
\cdot Z=2m$ since $C \cdot Z = 3m$ gives again a contradiction
from (\ref{eq*}). Then $0 \ge (C+2Z)^2= C^2 + 4Z^2 + 4 C \cdot Z
=-3m-8n+8m=5m-8n$; hence we get $m=n$. Let $x$ and $y$ be the two
intersection points and $C_x$(resp. $C_y$) the curve in $D$
containing $x$(resp. $y$). Then we have $D=C_x + 2C_y +D'$ with
$C_x \cdot C=n$ and $C_y \cdot C=n$ where $D'$ doesn't contain
$C_x$ nor $C_y$. Observe that we have $0=C_x \cdot Y=C_x^2+2C_x
\cdot C_y + C_x \cdot C + C_x \cdot D'=-2n+2C_x \cdot C_y +n+C_x
\cdot D'=-n+2C_x \cdot C_y+C_x \cdot D'$ ; hence $C_x \cdot C_y=0$
and $C_x \cdot D'=n$. We also have $0=C_y \cdot Y=C_y \cdot C_x +
2C_y^2 + C_y \cdot C + C_y \cdot D'=0-4n+n+C_y \cdot D'$ ; hence
$C_y \cdot D'=3n$. From $0=D' \cdot Y = D' \cdot C + D' \cdot C_x
+ 2 D' \cdot C_y + D'^2 = 0+n+6n+ D'^2$, we have $D'^2=-7n$. But
$D'^2 \equiv 0$ (mod $2n$).

(ii) Case $\beta=1$.

Let $x$ be the only intersection point. $D'$ denotes the sum of
curves in $D$ which don't meet $C$.

(Case 1) $D=C_{x,1} + C_{x,2} + C_{x,3}+ D'$ where $C_{x,i}$ meets
$C$ at $x$.

It is impossible since there is no cycle in $D$.

(Case 2) $D=C_{x,1} + C_{x,2} + D'$ with $C \cdot C_{x,1}=m$ and
$C \cdot C_{x,2}=2m$.

We have $m \ge n$ and $C \cdot Z=3m$ which contradict (\ref{eq*}).

(Case 3) $D=C_{x,1} + 2C_{x,2} + D'$ with $C \cdot C_{x,i}=m$.

We may assume that $C \cdot Z =2m$ since otherwise we get a
contradiction from ({\ref{eq*}}) along with $m \ge n$. Since $0
\ge (C+2Z)^2= C^2 + 4Z^2 + 4 C \cdot Z =-3m-8n+8m=5m-8n$, we have
$m=n$. Then $0=C_{x,1} \cdot Y=C_{x,1}^2+2C_{x,1} \cdot C_{x,2} +
C_{x,1} \cdot C + C_{x,1} \cdot D' \ge -2n+2C_{x,1} \cdot C_{x,2}
+ n \ge -2n+2n+n$ which is impossible.

(Case 4) $D=C_x + D'$ with $C \cdot C_x=3m$.

(a) $C_x^2=-2n$

We have $0=C_x \cdot Y= C_x \cdot C + C_x^2 + C_x \cdot D'\ge
3m-2n+n$. Since $n$ is a multiple of $3m$, we have $3m=n$ which is
impossible due to char($k$)=2.

(b) $C_x^2=-4n$

Similarly we have $3m=2n$. It is also impossible by the same
reason as in (a).

(Case 5) $D=3C_x + D'$ with $C \cdot C_x=m$.





Let me at first list all the possible configurations for $Y$ and
explain why we have them later.

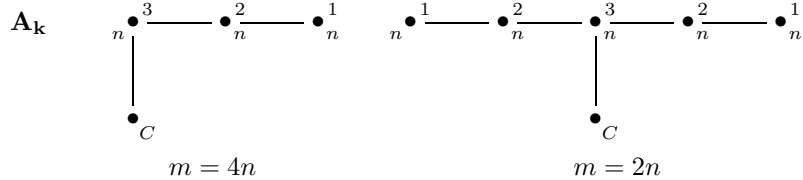
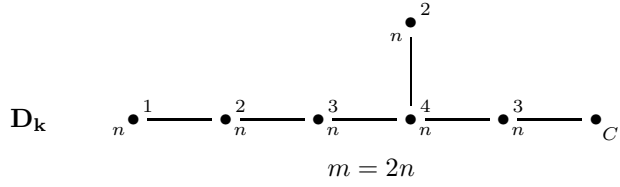
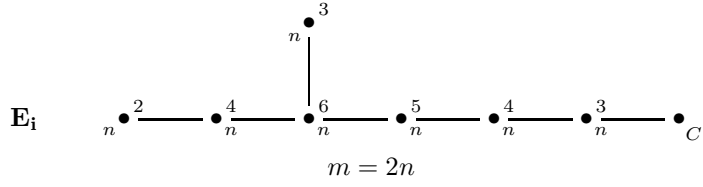
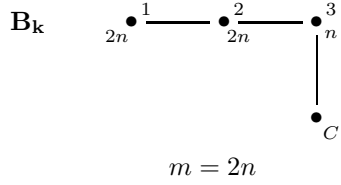
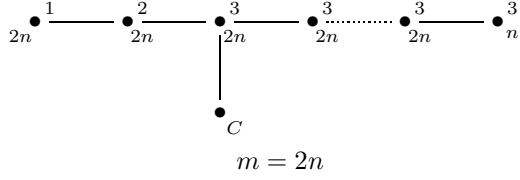
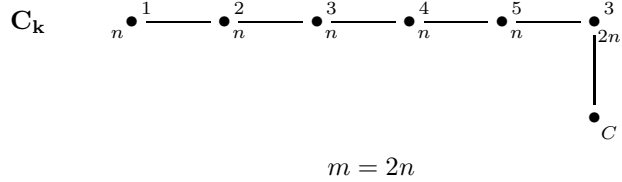
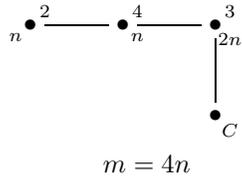
\xymatrix{ {\bf A_k} \ar@{}[r]_>{n} & {\bullet} \ar@{-}[r]^<{3}
\ar@{-}[d] & {\bullet} \ar@{-}[r]^<{2}_<{n} & {\bullet}
\ar@{}[r]^<{1}_<{n} \ar@{}[r]_>{n} & {\bullet} \ar@{-}[r]^<{1} &
{\bullet} \ar@{-}[r]^<{2}_<{n} & {\bullet} \ar@{-}[r]^<{3}_<{n}
\ar@{-}[d] & {\bullet} \ar@{-}[r]^<{2}_<{n} & {\bullet}
\ar@{}[r]^<{1}_<{n} &
\\
& {\bullet} \ar@{}[r]_<{C} & & & & & {\bullet} \ar@{}[r]_<{C} & }

\xymatrix{& & m=4n & & & &  m=2n &}

\vskip 3mm

\xymatrix{ & & & & {\bullet} \ar@{}[r]^<{2} & \\
{\bf D_k} \ar@{}[r]_>{n} & {\bullet} \ar@{-}[r]^<{1} & {\bullet}
\ar@{-}[r]^<{2}_<{n} & {\bullet} \ar@{-}[r]^<{3}_<{n} & {\bullet}
\ar@{-}[r]^<{4}_<{n} \ar@{-}[u]^>{n} & {\bullet}
\ar@{-}[r]^<{3}_<{n} & {\bullet} \ar@{}[r]_<{C} &  }

\xymatrix {& & & & m=2n}

\vskip 1.5mm

\xymatrix{ & & & {\bullet} \ar@{}[r]^<{3} & \\
{\bf E_i} \ar@{}[r]_>{n} & {\bullet} \ar@{-}[r]^<{2} & {\bullet}
\ar@{-}[r]^<{4}_<{n} & {\bullet} \ar@{-}[r]^<{6}_<{n}
\ar@{-}[u]^>{n} & {\bullet} \ar@{-}[r]^<{5}_<{n} & {\bullet}
\ar@{-}[r]^<{4}_<{n} & {\bullet} \ar@{-}[r]^<{3}_<{n} & {\bullet}
\ar@{}[r]_<{C} & }

\xymatrix {& & & & m=2n}

\vskip 3mm

\xymatrix{ {\bf B_k} \ar@{}[r]_>{2n} & {\bullet} \ar@{-}[r]^<{1} &
{\bullet} \ar@{-}[r]^<{2}_<{2n} &
{\bullet} \ar@{-}[d] \ar@{}[r]^<{3}_<{n} & \\
& & & {\bullet} \ar@{}[r]_<{C} & }

\xymatrix {& & m=2n &}

\vskip 1.5mm

\xymatrix { \hskip 4mm \ar@{}[r]_>{2n} & {\bullet} \ar@{-}[r]^<{1}
& {\bullet} \ar@{-}[r]^<{2}_<{2n} & {\bullet} \ar@{-}[d]
\ar@{-}[r]^<{3}_<{2n} & {\bullet} \ar@{.}[r]^<{3}_<{2n} &
{\bullet} \ar@{-}[r]^<{3}_<{2n} & {\bullet} \ar@{}[r]^<{3}_<{n} & \\
& & & {\bullet} \ar@{}[r]_<{C} & }

\xymatrix {& & & & m=2n}

\vskip 3mm

\xymatrix { {\bf C_k} \ar@{}[r]_>{n} & {\bullet} \ar@{-}[r]^<{1} &
{\bullet} \ar@{-}[r]^<{2}_<{n} & {\bullet} \ar@{-}[r]^<{3}_<{n} &
{\bullet} \ar@{-}[r]^<{4}_<{n} & {\bullet} \ar@{-}[r]^<{5}_<{n} &
{\bullet} \ar@{-}[d] \ar@{}[r]^<{3}_<{2n} & \\
& & & & & & {\bullet} \ar@{}[r]_<{C} & }

\xymatrix { & & & & m=2n}

\vskip 1.5mm

\xymatrix{& \ar@{}[r]_>{n} & {\bullet} \ar@{-}[r]^<{2} & {\bullet}
\ar@{-}[r]^<{4}_<{n} & {\bullet} \ar@{-}[d] \ar@{}[r]^<{3}_<{2n} & \\
& & & & {\bullet} \ar@{}[r]_<{C} & }

\xymatrix{& & & m=4n}

\vskip 3mm

(a) $C_x^2=-2n$

Claim. Any curve $C_i$ meeting with $C_x$ cannot have multiplicity
1.

{\it Proof.}  Suppose that $C_i$ has multiplicity 1 in $Y$. Then
we have $0=C_i \cdot Y \geq C_i^2+3C_i \cdot C_x$.  If
$[\Gamma(C_i,\oh_{C_i}):k]=n$(resp. $2n$), $C_i^2+3C_i \cdot
C_x=n$(resp. $2n$). It gives a contradiction. {\it (QED of Claim)}

We have $0=C_x \cdot Y=3C_x^2+C_x \cdot C+C_x \cdot D' \ge
-6n+m+n$ ; hence $m=n,2n$ or $4n$.

If $m=n$, we have $C_x \cdot D'=5n$ since $0=C_x \cdot Y=
3C_x^2+C_x \cdot C+C_x \cdot D'=-6n+m+ C_x \cdot D'$. Then we get
$D'^2=-15n$ from $D' \cdot Y=0$. It is impossible since $D'^2
\equiv 0$ (mod $2n$).

If $m=4n$, $C_x \cdot D' = 2n$; hence we have only one possibility
for $D$(note that $C_x$ must be an end since otherwise one of
curves meeting with $C_x$ has multiplicity 1).

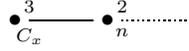
\xymatrix{  & {\bullet} \ar@{-}[r]_<{C_x} ^<{3} & {\bullet}
\ar@{.}[r]_<{n} ^<{2} & }


Let's use the condition that the intersection number of any curve
$C_i$ with $Y$ is zero. In the diagram the second curve should
meet only one curve,say $C_3$, with the multiplicity 1 and
$[\Gamma(C_3,\oh_{C_3}):k]=n$. It gives the first of $A_k$ in all
the possible configuration.

If $m=2n$, $C_x \cdot D' = 4n$; hence we have only the following
possibilities for $D$. Note that there is no case where $C_x$
meets three other curves in $D$. It is because of the above claim.

\xymatrix{ (1) & {\bullet} \ar@{-}[r]_<{C_x} ^<{3} & {\bullet}
\ar@{.}[r]_<{n} ^<{4} & }

\xymatrix{ (2) & {\bullet} \ar@{-}[r]_<{C_x} ^<{3} & {\bullet}
\ar@{.}[r]_<{2n} ^<{2} & }

\xymatrix{ (3) & \ar@{.}[r] & {\bullet} \ar@{-}[r]_<{n} ^<{2} &
{\bullet} \ar@{-}[r]_<{C_x} ^<{3} & {\bullet} \ar@{.}[r]_<{n}
^<{2} & }

Again use the intersection condition. Then we will have the $D_k$
case and the $E_i$ case from the diagram (1). Also we have the
first in $B_k$ from the diagram (2) and the second in $A_k$ from
the (3).

(b) $C_x^2=-4n$

We have $0=C_x \cdot Y=3C_x^2+C_x \cdot C+C_x \cdot D' \ge
-12n+m+2n$ ; hence $m=2n,4n$ or $8n$.

Claim. Any curve $C_i$ meeting with $C_x$ cannot have multiplicity
1 and if the multiplicity is 2, then
$[\Gamma(C_i,\oh_{C_i}):k]=2n$.

{\it Proof.}  Suppose that $C_i$ has multiplicity 1 in $Y$. Then
we have $0=C_i \cdot Y \geq C_i^2+3C_i \cdot C_x=C_i^2+6n$ which
is impossible. If $C_i$ has multiplicity 2 in $Y$, then $0=C_i
\cdot Y \geq 2C_i^2+3C_i \cdot C_x=2C_i^2+6n$. Therefore $C_i^2$
must be $-4n$. {\it (QED of Claim)}

If $m=8n$, $C_x \cdot D' = 4n$. It implies that $C_x$ must be an
end in $D$ since otherwise at least one curve meeting $C_x$ has
multiplicity 1. Hence we have only one case :

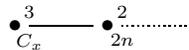
\xymatrix{ & {\bullet} \ar@{-}[r]_<{C_x} ^<{3} & {\bullet}
\ar@{.}[r]_<{2n} ^<{2} & }

The intersection property easily shows that the only possible
diagram for $D$ is

\xymatrix{ & {\bullet} \ar@{-}[r]_<{C_x}^<{3} & {\bullet}
\ar@{-}[r]_<{2n}^<{2} & {\bullet} \ar@{}[r]_<{2n}^<{1} & . } But
then we have a contradiction since $n$ is the smallest among
$n_i$'s in $D$.

If $m=4n$, $C_x \cdot D' = 8n$; hence we have the following cases
for $D$ :

\xymatrix{ (1) & {\bullet} \ar@{-}[r]_<{C_x} ^<{3} & {\bullet}
\ar@{.}[r]_<{n} ^<{4} & }

\xymatrix{ (2) & {\bullet} \ar@{-}[r]_<{C_x} ^<{3} & {\bullet}
\ar@{.}[r]_<{2n} ^<{4} & }

\xymatrix{ (3) & \ar@{.}[r] & {\bullet} \ar@{-}[r]_<{2n} ^<{2} &
{\bullet} \ar@{-}[r]_<{C_x} ^<{3} & {\bullet} \ar@{.}[r]_<{2n}
^<{2} & }

We get the second in $C_k$ from (1) and nothing from (2) \& (3).

If $m=2n$, $C_x \cdot D' = 10n$; hence we have the following cases
for $D$ :

\xymatrix{ (1) & {\bullet} \ar@{-}[r]_<{C_x} ^<{3} & {\bullet}
\ar@{.}[r]_<{n} ^<{5} & }

\xymatrix{ (2) & {\bullet} \ar@{-}[r]_<{C_x} ^<{3} & {\bullet}
\ar@{.}[r]_<{2n} ^<{5} & }

\xymatrix{ (3) & \ar@{.}[r] & {\bullet} \ar@{-}[r]_<{2n} ^<{2} &
{\bullet} \ar@{-}[r]_<{C_x} ^<{3} & {\bullet} \ar@{.}[r]_<{n}
^<{3} & }

\xymatrix{ (4) & \ar@{.}[r] & {\bullet} \ar@{-}[r]_<{2n} ^<{2} &
{\bullet} \ar@{-}[r]_<{C_x} ^<{3} & {\bullet} \ar@{.}[r]_<{2n}
^<{3} & }

We get the first in $C_k$ from (1), the second in $B_k$ from (4)
and nothing from (2) \& (3).

Let $\sigma \in H^0(\omega_{X/S}^{\otimes 2}|_Y)$ with $\sigma(x)
\neq 0$. Then $\sigma|_C \not\equiv 0$ and $Z(\sigma)=y_1+y_2$ or
$y$($y_1=y_2$ could happen).

Claim 1. There is a basis of $H^0(\omega_{X/S}^{\otimes 2}|_Y)$
consisting of $\{\sigma_1,...,\sigma_{\frac32m } \}$ such that
$\sigma_i(x) \neq 0$ and so these don't vanish on any component.

{\it Proof.}  We have $H^1(\omega_{X/S}^{\otimes 2}|_Y)=0$ by
Lemma \ref{van}. Then Riemann-Roch gives
$[H^0(\omega_{X/S}^{\otimes 2}|_Y):k]=\frac32m$. The existence of
such basis comes from the openness of the nonvanishing property in
the statement . {\it (QED of Claim 1)}

Claim 2. $[H^0(Y,\oh_Y):k] \leq n$

{\it Proof.}  Choose $C'$ in $D$ with $[\Gamma(C',\oh_{C'}):k]=n$.
Then we have an exact sequence $$0 \to H^0(Y-C',\oh_{Y-C'}(-C'))
\to H^0(Y,\oh_Y) \to H^0(C',\oh_{C'}) \to ...$$ Here
$H^0(Y-C',\oh_{Y-C'}(-C')) \cong H^1(Y-C',
\omega_{Y-C'}(C'))^{\ast} \cong H^1(Y-C',\omega_{X}(Y))^{\ast}$.
But $H^1(Y-C',\omega_{X}(Y))=0$ by Lemma \ref{van} because $C$ has
multiplicity 1 in $Y$; hence $B^2 < 0$ for any subcurve $B$ of
$Y-C'$. {\it (QED of Claim 2)}

We can assume that $Z(\sigma_j)=y_1+y_2$ for some $j$ (Suppose
not. Then the zero divisor of each $\sigma_i$ has the second type.
If $\sigma_{j_1}$ and $\sigma_{j_2}$ have
$Z(\sigma_{j_1})=Z(\sigma_{j_2})=y$, $\sigma_{j_1} / \sigma_{j_2}$
is in $H^0(Y,\oh_Y)$. Claim 2 implies that there are
$\sigma_{j_1}$ and $\sigma_{j_2}$ satisfying
$Z(\sigma_{j_i})=y_i(i=1,2)$ with $y_1 \neq y_2$. Then
$\{\sigma_{j_1},\sigma_{j_2}\}$ generate $\omega_{X/S}^{\otimes
2}$). Let $C'$ be a $(-2)$-curve in $Y$ with
$[\Gamma(C',\oh_{C'}):k]=n$. Consider an exact sequence
$$0 \longrightarrow \omega_{X/S}^{\otimes 2}(-Y+C')|_{C'} \longrightarrow
\omega_{X/S}^{\otimes 2}|_Y \longrightarrow \omega_{X/S}^{\otimes
2}|_{Y-C'} \longrightarrow 0.$$ This induces an exact sequence $0
\longrightarrow H^0(\omega_{X/S}^{\otimes 2}(-Y+C')|_{C'})
\longrightarrow H^0(\omega_{X/S}^{\otimes 2}|_Y) \longrightarrow
H^0(\omega_{X/S}^{\otimes 2}|_{Y-C'}) \longrightarrow
H^1(\omega_{X/S}^{\otimes 2}(-Y+C')|_{C'}) \longrightarrow 0$. We
have $\deg \omega_{X/S}^{\otimes 2}(-Y+C')|_{C'}) = C'^2$. It
implies that $H^0(\omega_{X/S}^{\otimes 2}(-Y+C')|_{C'})=0$ and
$H^1(\omega_{X/S}^{\otimes 2}(-Y+C')|_{C'}) \cong
\Gamma(C',\oh_{C'})$. Hence we have an exact sequence $$0
\longrightarrow H^0(\omega_{X/S}^{\otimes 2}|_Y) \longrightarrow
H^0(\omega_{X/S}^{\otimes 2}|_{Y-C'}) \longrightarrow
\Gamma(C',\oh_{C'}) \longrightarrow 0.$$ Therefore
$H^0(\omega_{X/S}^{\otimes 2}|_{Y-C'})$ has a basis
${\sigma_1,...,\sigma_{\frac32m},\tau_1,...\tau_n}$ where
$\sigma_j$'s extend to $Y$. But we have $H^1(\omega_{X/S}^{\otimes
2}|_{Y-C'}(-y_i))=0$ by Lemma \ref{van} hence
$H^0(\omega_{X/S}^{\otimes 2}|_{Y-C'}(-y_i))$ has dimension
$n+\frac{m}2$. Therefore at least one of $\sigma_j$'s doesn't
vanish at $y_i$. This concludes the proof. \epf

\begin{co}
$\omega_{X_0}^{\otimes m}$ is globally generated for $m \ge 2$.
\end{co}

\bpf Observe that $\omega_{X_0}^{\otimes m} \cong
\omega_X^{\otimes m}(m X_0)|_{X_0}$ and then we follow the same
proof as in Theorem \ref{main}. \epf

\vskip 10mm

\section{Canonical Rings and Koszul Cohomology}

Using Theorem \ref{main} together with an adapted K. Konno's
argument in \cite{konno}, we will show that $R(X_0,\omega_{X_0})=
\bigoplus_{m \geq 0} H^0(X_0,\omega_{X_0}^{\otimes m})$ is
generated by elements of degree up to 5. Konno's argument uses a
formal generalization of Green's Koszul cohomology.

\subsection{Koszul cohomology}(see \cite{konno} or \cite{green})

Let $D$ be a curve (i.e. an effective divisor) in a regular
surface $X$ over a fixed field $K$. Let $\mathcal{L}$ and
$\mathcal{M}$ be two line bundles on $D$ and $W$ a subspace of
$H^0(D,\mathcal{L})$. Then we have natural differentials(or Koszul
maps)
$$d_{i,j}:\wedge^i W \otimes H^0(D,\mathcal{M}+j \mathcal{L})
\to \wedge^{i-1} W \otimes H^0(D,\mathcal{M}+(j+1) \mathcal{L}).$$
Put $K_{i,j}(D,\mathcal{M},W)=Ker(d_{i,j})/Im(d_{i+1,j-1})$. For
convenience we put $K_{i,j}(D,\mathcal{M},\mathcal{L}) :=
K_{i,j}(D,\mathcal{M},H^0(D,\mathcal{L}))$.

\begin{lm}(Vanishing Theorem, 1.2.2 in \cite{konno}) \label{van2}
$K_{i,j}(D,\mathcal{M},W)=0$ in the following cases.

(1) $h^0(D,\mathcal{M}+j\mathcal{L})=0$.

(2) $h^0(D,\mathcal{M}+j\mathcal{L})=1$ and
$\mathcal{M}+j\mathcal{L} \cong \oh_D$.

\end{lm}

Assume that $W$ is base point free with $\dim_K W=r+1$. Then we
have an exact sequence $$0 \to \mathcal{E} \to W \otimes
\mathcal{O}_D \to \mathcal{L} \to 0,$$ where $\mathcal{E}$ is a
locally free sheaf of rank $r$ and $\wedge^r \mathcal{E} \simeq
\mathcal{O}_D(-\mathcal{L})$. Tensoring
$\mathcal{O}_D(\mathcal{M}+j\mathcal{L})$, we get $$0 \to E_{i,j}
\to \wedge^i W \otimes \mathcal{O}_D(\mathcal{M}+j\mathcal{L}) \to
E_{i-1,j+1} \to 0,$$ where $E_{i,j}=\wedge^i \mathcal{E} \otimes
\mathcal{O}_D(\mathcal{M}+j\mathcal{L})$. Then we have
$d_{i,j}=\wedge^i W \otimes H^0(D,\mathcal{M}+j \mathcal{L}) \to
H^0(D,E_{i-1,j+1}) \to \wedge^{i-1} W \otimes
H^0(D,\mathcal{M}+(j+1)\mathcal{L})$. By chasing diagrams we can
prove :

\begin{lm}(Duality, 1.2.1 in \cite{konno}) \label{dual}

$K_{i,j}(D,\mathcal{M},W)$ is dual to
$K_{r-1-i,2-j}(D,K_D-\mathcal{M},W)$.

\end{lm}








\begin{co}(1.2.3 in \cite{konno}) \label{k2}
Let $D$ be a curve on a surface which satisfies $H^0(D,-mK_D)=0$
for any positive integer $m$ and $\dim_K H^0(D,\oh_D)=1$. Put
$g'=\dim_K H^0(D,K_D)$. Assume that $2 K_D$ is generated by its
global sections. Then $K_{i,j}(D,\oh_D,2 K_D)^* \cong
K_{3g'-5-i,2-j}(D,K_D,2 K_D)$. Hence $K_{i,j}(D,K_D,2
K_D)=K_{i,j}(D,\oh_D,2 K_D)=0$ for $j \geq 3$ and $K_{i,2}(D,K_D,2
K_D)=0$ when $3g'-5-i \geq 1$.
\end{co}

\bpf The statement is clear from Lemma \ref{dual} along with Lemma
\ref{van2}. \epf

\begin{pr} \label{uptofour}
Let $D$ be a $1$-connected curve on a surface with $K_D$ nef and
$\dim_K H^0(D,K_D) \geq 2$. Assume that $\oh_D(mK_D)$ is globally
generated for $m \geq 2$. Then $R(D,K_D)=\oplus_{n \geq 0} H^0(D,n
K_D)$ is generated in degree up to $4$.
\end{pr}

\bpf We only need to show that $$H^0(D,2K_D)\otimes H^0(D,n K_D)
\to H^0(D,(n+2)K_D)$$ is surjective for $n \ge 3$. It follows from
Corollary \ref{k2}. \epf

Using the above Proposition and the result from the previous
section, we have the following :

\begin{pr}
Assume that the original settings for $X \to Spec(R)$ hold . If
$X_0$ is $1$-connected, then $R(X_0, \omega_{X_0})$ is generated
in degree up to $4$.
\end{pr}
\bpf Let $K=H^0(X_0,\oh_{X_0})$. Then $\dim_K H^0(D,K_D) \geq 2$
since $2g-2=(2\dim_K H^0(D,K_D)-2)[K:k]$ and $g \geq 2$.
Proposition \ref{uptofour} concludes the proof. \epf

We will extend this proposition to the general case where $X_0$ is
not necessarily $1$-connected. Then $X_0=rY$ for some integer $r$.

\begin{lm} \label{surj}
Let $m \ge 2$. Then $H^0(iY,\omega_{X_0}^{\otimes m}) \to
H^0(jY,\omega_{X_0}^{\otimes m})$ is surjective for any $i>j>0$.
\end{lm}

\bpf We need to show that $H^1((i-j)Y,\omega_{X_0}^{\otimes
m}(-jY))=0$. It follows from Lemma \ref{van}. \epf

\begin{tm}
Without the $1$-connected condition for $X_0$, $R(X_0,
\omega_{X_0})$ is generated in degree up to $5$.
\end{tm}

\bpf Let $m \geq 2$. Lemma \ref{surj} along with the exact
sequence $0 \to \oh_Y(\omega_{X_0}^{\otimes m}(-(i-1)Y)) \to
\oh_{iY}(\omega_{X_0}^{\otimes m}) \to
\oh_{(i-1)Y}(\omega_{X_0}^{\otimes m}) \to 0$ implies that, in
order to prove the surjectivity of $H^0(X_0,\omega_{X_0}^{\otimes
m}) \otimes H^0(X_0,\omega_{X_0}^{\otimes 2}) \to
H^0(X_0,\omega_{X_0}^{\otimes m+2})$, it suffices to see that
$H^0(Y,\omega_{X_0}^{\otimes m}(-lY)) \otimes
H^0(Y,\omega_{X_0}^{\otimes 2}) \to H^0(Y,\omega_{X_0}^{\otimes
m+2}(-lY))$ is surjective for all $0 \leq l \leq r-1$.

Claim 1. $K_{0,j}(Y,\omega_{X_0}^{\otimes
2}(-lY),\omega_{X_0}^{\otimes 2})=0$ for all $j \geq 2$.

{\it Proof.}  We have $K_{0,j}(Y,\omega_{X_0}^{\otimes
2}(-lY),\omega_{X_0}^{\otimes 2}) \cong K_{s-1,2-j}(Y,
\omega_{X_0}^{\otimes -1}((l+1)Y-X_0),\omega_{X_0}^{\otimes 2})^{\ast}$ by
Lemma \ref{dual} where $s+1$ is the dimension of $H^0(Y,\omega_{X_0}^{\otimes 2})$. The claim is true once
$H^0(Y,\omega_{X_0}^{\otimes 3-2j}((l+1)Y-X_0))=0$ by Lemma
\ref{van2}. But $H^0(Y,\omega_{X_0}^{\otimes 3-2j}((l+1)Y-X_0)) \cong
H^1(Y,\omega_{X_0}^{\otimes 2j-2}(-lY))^{\ast}=0$ for all $j \geq
2$ by Lemma \ref{van}. {\it (QED of Claim 1)}

Claim 2. $K_{0,j}(Y,\omega_{X_0}(-lY),\omega_{X_0}^{\otimes 2})=0$
for all $j \geq 3$.

{\it Proof.}  Similar argument as in Claim 1 leads to showing the
vanishing of $H^1(Y,\omega_{X_0}^{\otimes 2j-3}(-lY))$ for all $j
\geq 3$. It is true again by Lemma \ref{van}. {\it (QED of Claim 2
)}

Claim 1 and 2 imply that $$H^0(Y,\omega_{X_0}^{\otimes m}(-lY))
\otimes H^0(Y,\omega_{X_0}^{\otimes 2}) \to
H^0(Y,\omega_{X_0}^{\otimes m+2}(-lY))$$ is surjective for all $m
\geq 4$ and all $0 \leq l \leq r-1$. It ends the proof of this
theorem. \epf

\begin{co} \label{uptofive}
Let $f:X \to S$ be a relatively minimal fibration of genus $g \ge 2$. Then 
$\mathcal{R}(f)=\bigoplus_{m \geq 0} \mathcal{R}_m=\bigoplus_{m
\geq 0} f_*(\omega_{X/S}^{\otimes m})$ is generated by elements in
degree up to $5$.
\end{co}

\section{Remarks}

\subsection{}

By copying the statements in \cite{konno}, we can prove the
following :

Let $F$ be a fiber in a relatively minimal fibration of curves of
genus $g \ge 2$ over an algebraically closed field $k$. Then the
canonical ring $R(F,K_F)$ is generated in degrees 1,2 and 3 except
when $F$ is a multiple fiber which contains a $(-1)$-elliptic
cycle $E$ such that $E \subset Bs|K_F|$. In the exceptional case,
it needs one more generator in degree 4.

In general, the methods in \cite{konno} cannot be applied ; e.g.
the third part of Lemma 1.2.2. in \cite{konno} doesn't hold if the
base field is not algebraically closed.

\subsection{}

I'd like to see a generalization of Corollary \ref{uptofive} to
the case where $S$ is not necessarily one dimensional. For this we
will need a similar global generatedness of $\omega_{X/S}^{\otimes
m}$. Here is one notable thing with $\dim S=2$. If we have a
family $f : X \to S$ of curves over a surface over an
algebraically closed field $k$ with $X$ and $S$ regular,
$\omega_{X/S}$ can be shown to be relatively nef using deformation
theory (See e.g. \cite{kollar}). Assuming $S$ to be projective
over $k$ and $A$ is sufficiently ample on $S$,
$\omega_{X/S}^{\otimes 2} \otimes f^{\ast}(A)$ is nef and big.
Using Keel's results in \cite{keel}, $\omega_{X/S}^{\otimes 2}
\otimes f^{\ast}(A)$ is semi-ample if $k$ is the algebraic closure
of a finite field.

\subsection{}

By Corollary \ref{uptofive} we can define $Proj(\mathcal{R}(f))$.
The map $X \to Proj(\mathcal{R}(f))$ is essentially the
contraction of all $(-2)$-chains in each fiber to points. There is
a very similar phenomenon for a regular minimal surface over
$\mathbb{C}$ (See \cite{bpv}).

\end{document}